\documentclass[12pt]{iopart}
\usepackage{iopams}
\usepackage{amsthm}
\usepackage[latin1]{inputenc} 
\usepackage{bbm}
\usepackage{cite}
\usepackage{amsfonts,amssymb}

\newcommand{\alg}{\mathcal{F}}
\newcommand{\prob}{\mathbb{P}}
\newcommand{\att}{\mathbb{E}}

\newcommand{\Xcont}{\{X_t\}_{t\ge0}}

\newcommand{\EM}{{}^q\!M^\tau_N}
\newcommand{\EMk}{{}^q\!M^\tau_{N_k}}

\newcommand{\EMone}{{}^q\!M^1_N}
\newcommand{\EMonek}{{}^q\!M^1_{N_k}}
\newcommand{\lcm}{\mbox{lcm}(\mathcal{T})}

\newtheorem{theorem}{Theorem}
\newtheorem{proposition}{Proposition}
\newtheorem{lemma}{Lemma}

\begin{document}

\title{Apparent multifractality of self--similar L\'evy processes}

\author{Marco Zamparo} 

\address{Dipartimento Scienza Applicata e Tecnologia and CNISM, Politecnico di
Torino, Corso Duca degli Abruzzi 24, I--10129 Torino, Italy} 


\ead{marco.zamparo@polito.it}

\begin{abstract}
Scaling properties of time series are usually studied in terms of the
scaling laws of empirical moments, which are the time average
estimates of moments of the dynamic variable. Nonlinearities in the
scaling function of empirical moments are generally regarded as a sign
of multifractality in the data. We show that, except for the Brownian
motion, this method fails to disclose the correct monofractal nature
of self--similar L\'evy processes. We prove that for this class of
processes it produces apparent multifractality characterised by a
piecewise--linear scaling function with two different regimes, which
match at the stability index of the considered process.  This result
is motivated by previous numerical evidence. It is obtained by
introducing an appropriate stochastic normalisation which is able to
cure empirical moments, without hiding their dependence on time, when
moments they aim at estimating do not exist.
\end{abstract}


\maketitle
\section{Introduction}

Many variables in natural and economic sciences appear to have some
scale invariance properties.  Since the beginning of fractal geometry
with the work of Mandelbrot in 1960's, statistical self--similarity
and scaling laws have been discovered in turbulence \cite{BPPV1984},
ecology \cite{HKG1999}, hydrology \cite{GW1993}, network data traffic
\cite{LTWW1994}, and finance \cite{DM2007} among other fields.  In the
specific context of time series, time--scale invariance for the
moments of the dynamic variable is generally observed \cite{DM2007}.
The scaling law of different moments is completely described by a
single exponent for monofractal dynamics, such as those corresponding
to self--similar processes, whereas a continuous spectrum of exponents
is needed for multifractal evolutions.  The monofractal and the
multifractal scaling properties are commonly regarded as fundamentally
different, the former arising from additive and the latter from
multiplicative processes \cite{MFC1997,CF2002}.

Although multifractality is very attractive, there is a certain
controversy over its emergence and nature, especially in
finance. While many authors have provided evidence of multifractality
in stock markets using different methods on the one hand
\cite{SSL1999,XG2003,WH2005,JZ12008,ZFTPGR2009,GHH2014}, several
papers report on the other hand that no sign of multifractal scaling
laws is found in financial time series \cite{JZ22008,Z2009,Z2012}. In
addition to this, spurious multifractal effects have been documented
for different model types that in principle do not account for such a
form of scaling properties
\cite{CG2000,BPM2000,KCH2006,HS2008,H2009,N12010,N22010,GNR2012,GL2014}.
This casts doubt on whether the observed multifractality always
corresponds to a genuine multifractal process.

Multifractality in time series is detected by inspecting the power law
in time of empirical moments, which are the time average estimates of
moments, and evaluating their scaling function for various orders
\cite{CF2002}. A linear scaling function corresponds to a monofractal
process, whereas a nonlinear scaling function should denote the
occurrence of multifractality.  The basic assumption underlying this
protocol is that moments exist. If this assumption is violated, then
unforeseen and misleading results may be obtained, such as in the case
of self--similar L\'evy processes different from the Brownian motion.
In fact, for these processes the scaling function of empirical moments
results to be piecewise linear with two different regimes when the
range of existing moments is supposed to be not known in advance, as
argued in \cite{SSL1999} and \cite{CG2000} on the basis of numerical
simulations. A similar accident occurs for the fractional L\'evy
motion \cite{HS2008}.  Piecewise linearity of the scaling function
would suggest a multifractal model even though self--similar L\'evy
processes and the fractional L\'evy motion are exactly monofractal
processes.  The nonlinearity of the scaling function is entirely
caused by the non--existence of some moments, which for both model
classes is due to the presence of fat tails in the distribution of the
dynamic variable.  This suggests that empirical facts considered as
hallmarks of multifractal systems may be reproduced by plain
fat--tailed processes, so that an apparent multifractality may be
consistent with the latter.  This circumstance deserves attention in
finance, where variables are widely believed to possess fat tails
\cite{C2001} and where scaling functions approximately exhibiting
piecewise linearity with two different regimes have already been found
in several daily exchange rates \cite{MGFRS2004}. Interestingly,
scaling functions which share this same shape have also been found in
the context of the natural sciences, such as solar wind data
\cite{WCHCFG2005}, diffusion in living cells \cite{GW2010}, and
transport phenomena in optical lattices \cite{DL2012}.

Discerning real from apparent multifractality demands a thorough study
of the estimators of moments with the purpose of isolating the
mechanisms that may lead to spurious effects. The probability
distribution of these estimators has been obtained for self--similar
L\'evy processes and the fractional L\'evy motion in the limit of long
time series, where scaling laws in time of the distribution emerge
\cite{HS2008}.  Within a special protocol, where the time horizon to
which scaling laws should refer is allowed to grow as the sample size
increases, similar results have been achieved for processes with more
general independent and identically distributed fat--tailed increments
\cite{H2009}.  The latter work has very recently been extended to
models with stationary and weakly dependent increments \cite{GL2014}.
These studies are however only partially informative as one must keep
an eye on the fact that a deep insight into apparent multifractality
of time series requires to understand whether or not and how the
single empirical moments can exhibit multifractal scaling laws when
the underlying process is not a multifractal process, as observed in
\cite{SSL1999} and \cite{CG2000}.  The scaling properties of the
single empirical moments are not directly related to those of their
distribution, so that invoking the latter does not help to unravel the
knot.  To this aim, stronger results are needed and should be reached
by somehow improving the weak convergence in distribution to the
almost sure convergence, which is the only one that can explain what
occurs in one single experiment.

In this paper we reconsider the problem of apparent multifractality of
self--similar L\'evy processes. We show that, except for the Brownian
motion, empirical moments systematically exhibit a multifractal
scaling law for almost all the sample paths of the underlying
self--similar L\'evy process, mistakenly leading to consider it a
multifractal model.  While non--degenerate distributions and their
scaling properties can be obtained by resorting to norming constants
only, the strategy we use to highlight the scaling laws of the single
empirical moments when moments do not exist is employing a stochastic
normalisation.  A stochastic normalisation is made necessary by an
important result due to Feller on sums of independent and identically
distributed random variables that are not integrable.  This strategy
allows us to completely explain the numerical findings of
\cite{SSL1999} and \cite{CG2000} by demonstrating that empirical
moments of any order acquire their own deterministic scaling laws with
a well--defined scaling function in the limit of long time series.
The scaling function is piecewise--linear with two different regimes
matching at the stability index of the underlying process.

The paper is organised as follows.  In Section \ref{Scalingprop} we
recall some fundamental notions that are concerned with multifractal
processes and the statistical methods for estimating their scaling
function. In Section \ref{Levyflight} we introduce self--similar
L\'evy processes and we review previously known facts about their
apparent multifractality.  We supplement the literature showing
impossibility of disclosing multifractal scaling laws of empirical
moments by resorting to norming constants only. Then we discuss a
natural stochastic candidate for normalisation purposes and we state
our main result on the apparent multifractality of self--similar
L\'evy processes.  The proof of this result is outlined afterwards,
postponing the most technical details that are needed in the
appendices in order not to interrupt the flow of the
presentation. Conclusions and prospects for future research are
finally reported in Section \ref{conclusion}.

\section{Scaling properties of stochastic processes and empirical analysis}
\label{Scalingprop}

Let on the probability space $(\Omega,\alg,\prob)$\footnote{As usual,
  $\Omega$ denotes the sample space, $\alg$ is a $\sigma$-algebra of
  events on $\Omega$, and $\prob$ is a probability measure on
  $\alg$. A property holds {\it almost surely} (a.s.\ for short) if it
  holds for all the samples $\omega$ belonging to an event
  $E\in\mathcal{F}$ with $\prob[E]=1$.} be given a continuous--time
real--valued stochastic process $\Xcont$ satisfying $X_0=0$ and having
stationary increments. We recall that the {\it increment} of $\Xcont$
at time $t$ and scale $s$, or over the time window from $t$ to $t+s$,
is the variable $\Delta_t^sX:=X_{t+s}-X_t$ and that the process
possesses {\it stationary increments} if
$\Delta_t^sX\stackrel{d}{=}\Delta_0^sX$ for all $t>0$ and $s>0$, where
equality is in distribution.  Scale--invariance properties of $\Xcont$
can be defined by specifying the scaling rules of moments
\cite{MFC1997,CF2002}.  The process is said to be {\it multifractal}
if there exist some real numbers $T>0$ and $b>0$ and some functions
$\mu$ and $\nu$ with domain $[0,b)$ such that the scaling law
\begin{equation}
\att\big[|X_t|^q\big]=\mu(q) \cdot t^{\nu(q)}
\label{multifractal}
\end{equation}
holds for all $t\le T$ and non--negative $q<b$.  The quantity
$\att[|X_t|^q]$ is the {\it moment} of the process at time $t$ and
order $q$.  The function $\nu$ is the {\it scaling function} and
necessarily $\nu(0)=0$ and $\nu$ is concave \cite{MFC1997,CF2002}. The
simplest multifractal processes are characterised by a linear scaling
function and are also termed {\it monofractal}.  Self--similar
processes constitute an example of monofractal processes, whereas an
example of multifractal processes showing non linearities in the
scaling function is represented by the so--called log--infinitely
divisible multifractal processes \cite{MB2002,BM2003}.  We recall that
the process $\Xcont$ is said to be {\it self--similar} if there exists
a {\it scaling exponent} $H>0$ such that $X_{at}\stackrel{d}{=}a^H
X_t$ for all $a>0$ and $t>0$.  Self--similar processes verify
(\ref{multifractal}) for all $t$ with $b=\sup\{q\ge
0:\att[|X_1|^q]<\infty\}$, $\mu(q)=\att[|X_1|^q]$, and $\nu(q)=Hq$.

Empirical scaling analysis aims at assessing scaling properties of the
process $\Xcont$ on the basis of time series data and requires
estimation of the moment $\att[|X_\tau|^q]$ for several times $\tau$
and orders $q$.  The value $\nu(q)$ of the scaling function at $q$ is
evaluated by linear regression of logarithms of estimated moments on
$\ln \tau$ \cite{CF2002}.  In the typical experiment, $N+1$
measurements of the process are recorded at the equally spaced time
instants $0=\tau_0<\tau_1<\cdots<\tau_N$ and the times $\tau$
identifying moments to be estimated are multiples of the lag between
consecutive measurements. We choose unit of time in such a way that
$\tau_1=1$, so that $\tau_n=n$ for all $n$ and times $\tau$ are
integers.  Within this framework, the {\it estimator} $\EM$ of the
moment $\att[|X_\tau|^q]$ is defined as the {\it sample average} of
the increments of $\Xcont$ over $N/\tau$ non--overlapping, consecutive
time windows of size $\tau$:
\begin{equation}
\EM:=\frac{1}{N/\tau}\sum_{n=0}^{N/\tau-1}\bigl|\Delta_{\tau n}^\tau X\bigr|^q=\frac{1}{N/\tau}\sum_{n=1}^{N/\tau}\bigl|X_{n\tau}-X_{(n-1)\tau}\bigr|^q.
\label{Mestim}
\end{equation}
The realisation $\EM(\omega)$ of the random variable $\EM$
corresponding to the available data
$0=X_0(\omega),X_1(\omega),\ldots,X_N(\omega)$ is called {\it
  empirical moment} and precisely constitutes the {\it estimate} of
$\att[|X_\tau|^q]$. With the purpose of estimating the moment
$\att[|X_\tau|^q]$ for all $\tau$ up to a certain given integer
temporal horizon $\mathcal{T}$, the number of measurements $N$ is
assumed to be a multiple of the {\it least common multiple} $\lcm$ of
the first $\mathcal{T}$ integers. This way, the number $N/\tau$ of
consecutive increments that appears in (\ref{Mestim}) is an integer
for all $\tau\le\mathcal{T}$. In practical situations, $N/\mathcal{T}$
has to be large in order to collect enough statistics.

The hypothesis of stationary increments yields that
$\att[\EM]=\att[|X_\tau|^q]$ for all those times $\tau\le\mathcal{T}$
and orders $q$ that satisfy $\att[|X_\tau|^q]<\infty$, so that $\EM$
is an {\it unbiased estimator}. We notice that\footnote{As
  $X_{\tau+1}=\sum_{n=0}^\tau\Delta_n^1X$ and $\bigl|\sum_{n=0}^\tau
  \delta_n\bigr|^q\le 2^{q\tau}\cdot\sum_{n=0}^\tau|\delta_n|^q$ for
  all $q\ge 0$ and $\delta_0,\ldots,\delta_\tau$, the equality in law
  $\Delta_n^1X\stackrel{d}{=}X_1$ gives that $\att[|X_\tau|^q]<\infty$
  for each integer $\tau$ if and only if $\att[|X_1|^q]<\infty$.}
$\att[|X_\tau|^q]<\infty$ for all integers $\tau$ if and only if
$\att[|X_1|^q]<\infty$. The estimation procedure is said to be {\it
  consistent} if for all $q$ such that $\att[|X_1|^q]<\infty$
\begin{equation*}
\prob\Bigg[\bigcap_{\tau=1}^{\mathcal{T}}\biggl\{\omega\in\Omega:\lim_{k\uparrow\infty}\EMk(\omega)=\att\bigl[|X_\tau|^q\bigr]\biggr\}\Bigg]=1,
\end{equation*}
where $N_k:=k\,\lcm$ is a multiple of $\lcm$.  Consistency of
estimation is the only property that supports the use of empirical
moments and means, in a nutshell, that empirical moments converges to
moments as the number of data points increases indefinitely. We stress
that moments estimation requires that the range of existing moments is
known. Otherwise, one runs the risk of evaluating scaling functions
beyond the range of validity of moments scaling laws obtaining
possible misleading results, such as spurious multifractality. In this
respect, one must be aware that the empirical moment $\EM(\omega)$ is
well--defined for all $q\ge 0$ and does not allow to distinguish those
$q$ for which $\att[|X_\tau|^q]<\infty$ from those for which
$\att[|X_\tau|^q]=\infty$.

\section{The case of self--similar L\'evy processes}
\label{Levyflight}

A {\it L\'evy process} $\Xcont$ is a stochastic process with
stationary independent increments such that $X_0=0$. {\it Independent
  increments} means that the increments of the process over
non--overlapping time windows are independent variables. It is
well--known that a L\'evy process $\Xcont$ is self--similar if and
only if the distribution of $X_t$ is strictly stable \cite{Sato}. {\it
  Strictly stability} is expressed in terms of characteristic
functions as
\begin{equation*}
\att\bigl[\exp(i k X_t)\bigr]=\exp\bigl[-t\psi(k)\bigr],
\end{equation*}
where, with admissible parameters $0<\alpha\le 2$, $\sigma>0$, and
$-1\le\gamma\le 1$ when $\alpha\ne 1$ and any $\gamma$ if $\alpha=1$,
the function $\psi$ is defined as\footnote{As usual, $\mbox{sgn}(k)$
  denotes the sign of $k$: $\mbox{sgn}(k):=1$ if $k\ge 0$ and
  $\mbox{sgn}(k):=-1$ if $k<0$.}
\begin{equation*}
\psi(k):=\cases{
\sigma^\alpha|k|^\alpha\Bigl[1-i\gamma\tan\Bigl(\frac{\pi\alpha}{2}\Bigr)\mbox{sgn}(k)\Bigr] & if $\alpha\ne 1$;\\
\sigma|k|-i\gamma k & if $\alpha=1$.\\
}
\end{equation*}
The parameter $\alpha$ is called the {\it stability index} and
determines the scaling exponent $H$ associated to a self--similar
L\'evy process as $H=1/\alpha$. When $\alpha=2$ the distribution of
$X_t$ reduces to a Gaussian distribution and the process $\Xcont$
becomes a version of the Brownian motion. Parameters $\sigma$ and
$\gamma$ are a {\it scale parameter} and a {\it skewness parameter},
respectively. If $\gamma=0$, then $-X_t$ is distributed as $X_t$ for
each $t\ge 0$ and the self--similar L\'evy process is said to be {\it
  symmetric}.

\subsection{Main properties}

For the purposes of the present work, the most important feature of a
self--similar L\'evy process $\Xcont$ with stability index $\alpha$ is
that the distribution of $X_t$ displays {\it fat tails} with {\it tail
  index} $\alpha$ when $\alpha<2$ and $t>0$. Hereafter we implicitly
assume that $\alpha<2$. A power--law decay with exponent $\alpha$ of
both the right tail distribution $\prob[X_t>x]$ and the left tail
distribution $\prob[X_t<-x]$ can be immediately found out in the
simplest case $\alpha=1$, since in this case\footnote{If $f$ and $g$
  are real functions of a real variable $x$ defined for all large
  enough $x$, $f(x)\sim g(x)$ means that
  $\lim_{x\uparrow\infty}f(x)/g(x)=1$.}
\begin{equation*}
\prob[X_t>x]=\prob[X_t<-x]=\frac{1}{2}-\frac{1}{\pi}\arctan\biggl(\frac{x-t\gamma}{t\sigma}\biggr)\sim\frac{2t}{\pi}\frac{\sigma}{x}.
\end{equation*}
The case $\alpha\ne 1$ is a little bit more involved as the right tail
distribution shows a fat tail with tail index $\alpha$ only if
$-1<\gamma\le 1$ (see \cite{Zolotarev}, Theorem 2.4.2 if $\alpha<1$
and Corollary 2 of Theorem 2.5.1 if $\alpha>1$). On the contrary, if
$\gamma=-1$, then $\prob[X_t>x]=0$ for all $x>0$ when $\alpha<1$ or
$\prob[X_t>x]$ shrinks exponentially with $x>0$ when $\alpha>1$ (see
\cite{Zolotarev}, Theorem 2.5.2). For $\gamma>-1$ we exactly find
\begin{equation*}
\prob[X_t>x]\sim\frac{t}{\pi}(1+\gamma)\Gamma(\alpha)\sin\Bigl(\frac{\pi\alpha}{2}\Bigr)\Bigl(\frac{\sigma}{x}\Bigr)^{\alpha},
\end{equation*}
where $\Gamma$ is the Euler's gamma function.  Symmetrically, the left
tail distribution decays as a power law with exponent $\alpha$ if
$-1\le\gamma<1$, whereas $\prob[X_t<-x]=0$ for each $x>0$ or
$\prob[X_t<-x]$ decreases exponentially with $x>0$ according as
$\alpha<1$ or $\alpha>1$ if $\gamma=1$. When $\gamma<1$ we have that
\begin{equation*}
\prob[X_t<-x]\sim\frac{t}{\pi}(1-\gamma)\Gamma(\alpha)\sin\Bigl(\frac{\pi\alpha}{2}\Bigr)\Bigl(\frac{\sigma}{x}\Bigr)^{\alpha}.
\end{equation*}
In spite of the complexity introduced by the skewness parameter
$\gamma$, combining all together these results we recognise that the
plain power law
\begin{equation}
\prob[|X_t|>x]=\prob[X_t>x]+\prob[X_t<-x]\sim \frac{tc}{x^\alpha}
\label{tailabs}
\end{equation}
holds for all $\alpha<2$ and admissible value of other parameters with
\begin{equation}
c:=\frac{2}{\pi}\Gamma(\alpha)\sin\biggl(\frac{\pi\alpha}{2}\biggr)\sigma^{\alpha}>0.
\label{c_def}
\end{equation}
The consequence of (\ref{tailabs}) is that $\att[|X_t|^q]<\infty$ for
each $t>0$ if and only if $q<\alpha$.

Two other relevant features of the self--similar L\'evy process
$\Xcont$ stem from the fact that the sequence $\{\Delta_{\tau n}^\tau
X\}_{n\ge 0}$ of the increments over non--overlapping, consecutive
time windows of size $\tau$ is a sequence of independent and
identically distributed (i.i.d.\ for short) random variables. A first
consequence of this fact is that the {\it strong law of large numbers}
applies ensuring that for each integer temporal horizon $\mathcal{T}$
and non--negative order $q<\alpha$
\begin{equation}
\prob\Biggl[\bigcap_{\tau=1}^{\mathcal{T}}\biggl\{\omega\in\Omega:\lim_{k\uparrow\infty}\EMk(\omega)=\att\bigl[|X_\tau|^q\bigr]=\mu(q)\cdot\tau^{\nu(q)}\biggr\}\Biggr]=1,
\label{conv_emp_moment}
\end{equation}
where $\mu(q)=\att[|X_1|^q]$ and $\nu(q)=q/\alpha$. This means that
the estimation scheme previously discussed is consistent on the one
hand, and that empirical moments inherit the monofractal nature of
moments for orders smaller than the stability index on the other hand.
Another consequence of the fact that $\{\Delta_{\tau n}^\tau X\}_{n\ge
  0}$ is an i.i.d.\ sequence is a simple equality in law relating the
estimators $\EM$ and $\EMone$. Indeed, the self--similarity of the
process $\Xcont$ yields $\Delta_{\tau n}^\tau
X\stackrel{d}{=}\Delta_0^\tau
X\stackrel{d}{=}X_\tau\stackrel{d}{=}\tau^{1/\alpha}X_1
\stackrel{d}{=}\tau^{1/\alpha}\Delta_n^1X$ for all $n\ge 0$. Then,
bearing in mind the definition (\ref{Mestim}) of the estimator $\EM$,
we find for each non negative $q$, integer $\tau\le \mathcal{T}$, and
$N$ multiple of $\lcm$ the equality in law
\begin{equation}
\EM~\stackrel{d}{=}~\frac{\tau^{q/\alpha}}{N/\tau}\sum_{n=0}^{N/\tau-1}\bigl|\Delta_n^1X\bigr|^q=\tau^{q/\alpha}\cdot{}^q\!M^1_{N/\tau}.
\label{dist_Mt_M1}
\end{equation}

\subsection{Evidence of apparent multifractality}
\label{Evidence_mul}

Focusing on symmetric self--similar L\'evy processes, in \cite{CG2000}
the authors probed the empirical analysis to evaluate the scaling
function by means of numerical simulations. Mimicking practical
situations where no {\it a priori} information is available, they
pretended not to know the stability index $\alpha$ and attempted to
estimate moments on synthetic time series for trial orders $q$ up to
and above the threshold $\alpha$. Such an attempt is always possible
as the estimator $\EM$ is well--defined for each non negative $q$ but
estimation is phony for $q\ge\alpha$ as $\att[|X_\tau|^q]=\infty$ for
all $\tau>0$ in this case.  Surprisingly, the authors found that
empirical moments of any order obey an approximate scaling law in time
with a resulting scaling function which is concave and independent of
the sample.  They argued that such a scaling property becomes exact in
the limit of long time series and extrapolated the piecewise--linear
{\it empirical scaling function}
\begin{equation}
\nu_e(q):=\cases{
q/\alpha & if $q<\alpha$; \\
1 & if $q\ge\alpha$.\\
}
\label{empiricalnu}
\end{equation}
The same observation was previously reported in \cite{SSL1999}.  While
scaling laws at orders $q<\alpha$ were expected in view of
(\ref{conv_emp_moment}), the unforeseen behaviour of empirical moments
in the regime $q\ge\alpha$ was somehow ascribed to finite sample
effects.  The important lesson to be learnt from this study was that
if one trusts that empirical moments correctly estimate moments of
given trial orders for time series generated by self--similar L\'evy
processes of unknown stability index, then one is led to think that
data correspond to a non--trivial multifractal process even though
they do not. {\it Apparent multifractality} can then emerge when
fat--tailed distributions are associated to data.

Understanding apparent multifractality of self--similar L\'evy
processes from a mathematical standpoint amounts to explain whether or
not and why empirical moments exhibit scaling laws in time at all
orders for almost all the time series generated by the process.  As a
first speculative contribution to the problem, the distribution of the
estimator $\EM$ was investigated in \cite{HS2008} for $q>\alpha$.  The
idea was resorting to the {\it norming constant} $a_N:=N^{q/\alpha-1}$
in order to get in the large $N$ limit at a non--degenerate
distribution for $\EM/a_N$ whose dependence on $\tau$ may be disclosed.
The following theorem reports and completes with the case $q=\alpha$
the convergence in distribution presented in \cite{HS2008}. The proof
is provided in the next paragraph.
\begin{theorem}
\label{th:conv}
Let $\Xcont$ be a a self--similar L\'evy process of stability index
$\alpha<2$ defined on some probability space
$(\Omega,\mathcal{F},\prob)$ and let $c$ be the constant defined by
(\ref{c_def}). Let $\tau\le \mathcal{T}$ be positive integers and let
{\upshape $N_k:=k\,\lcm$} be a multiple of {\upshape $\lcm$}. Then
\begin{equation*}
\EMk -c\,\tau\ln\biggl(\frac{N_k}{\tau}\biggr)\stackrel{d}{\longrightarrow}\tau\cdot Z_q
\end{equation*}
as $k\uparrow\infty$ if $q=\alpha$ and
\begin{equation*}
N_k^{1-q/\alpha}\cdot\EMk\stackrel{d}{\longrightarrow}\tau\cdot Z_q
\end{equation*}
as $k\uparrow\infty$ if $q>\alpha$, $Z_q$ being a non--degenerate
stable variable of stability index $\alpha/q$ with law independent of
$\tau$.
\end{theorem}

Theorem \ref{th:conv} shows that the distribution of $\EM/a_N$ is
endowed with well--defined scaling laws\footnote{Similar results were
  also obtained for the fractional L\'evy motion \cite{HS2008} and for
  processes with more general fat--tailed independent \cite{H2009} and
  weakly dependent \cite{GL2014} increments within a special
  estimation scheme where, for a given $\xi\in(0,1)$, $\tau$ grows as
  $N^\xi$ as $N$ increases. In \cite{H2009} and \cite{GL2014}
  convergence in distribution was improved to convergence in
  probability.} in the size $\tau$ of time windows for orders $q$
strictly larger than $\alpha$, where no centering procedure is
required to reach a non degenerate limit when $N$ is sent to infinite.
Interestingly, the exponent involved in these laws coincides with that
found in \cite{SSL1999} and \cite{CG2000}.  However, despite the nice
result, this theorem does not give reasoning for the behaviour of
empirical moments observed in \cite{CG2000} since it only deals with
their statistics over the different samples on the one hand and since
it is not obvious that their own scaling in time may be deduced from
scaling laws of their distribution on the other hand.  Actually, we
believe that the basic idea of employing a norming constant is not
sufficient to explain apparent multifractality of self--similar L\'evy
processes.  The reason is that, contrary to the case of distributions,
when $q\ge\alpha$ no norming constant $a_N$ exists with the property
that for almost all $\omega\in\Omega$ the normed empirical moment
$\EM(\omega)/a_N$ attains at large $N$ a somehow useful limit to
detect scaling laws. This fact is made evident by the following
proposition, which is proven in the next paragraph.
\begin{proposition}
\label{no_cost}
Assume that $\Xcont$ is a self--similar L\'evy process of stability
index $\alpha<2$ defined on some probability space
$(\Omega,\mathcal{F},\prob)$. Let $\mathcal{T}$ be a positive integer
and let {\upshape $N_k:=k\,\lcm$} be a multiple of {\upshape
  $\lcm$}. Fix any $q\ge\alpha$ and let $\{a_{N_k}\}_{k\ge 1}$ be an
increasing sequence of positive numbers.  Then
\begin{equation*}
\prob\Biggl[\bigcap_{\tau=1}^{\mathcal{T}}\biggl\{\omega\in\Omega:\lim_{k\uparrow\infty}\frac{\EMk(\omega)}{a_{N_k}}=0\biggr\}\Biggr]=1
\end{equation*}
or
\begin{equation*}
\prob\Biggl[\bigcap_{\tau=1}^{\mathcal{T}}\biggl\{\omega\in\Omega:\limsup_{k\uparrow\infty}\frac{\EMk(\omega)}{a_{N_k}}=\infty\biggr\}\Biggr]=1
\end{equation*}
according as the series $\sum_{k=1}^\infty (ka_{N_k})^{-\alpha/q}$
converges or diverges.
\end{proposition}

Proposition \ref{no_cost} states that, whatever $a_N$ is, if
$q\ge\alpha$ either $\EM/a_N$ converges a.s.\ to zero as $N$
is sent to infinite or it does not converge at all.  As a consequence,
no property of single empirical moments can be highlighted by
resorting to a norming constant $a_N$ when $q\ge\alpha$, not even when
$a_N$ coincides with $N^{q/\alpha-1}$ as in Theorem \ref{th:conv}.

\subsection{Proof of Theorem \ref{th:conv} and Proposition \ref{no_cost}}
\label{proof:theorem_and_proposition}

In order to prove Theorem \ref{th:conv} it is enough to show that
there exists a non--degenerate stable variable $Z_q$ of stability
index $\alpha/q$ such that
\begin{equation}
\EMone-c\ln(N)\stackrel{d}{\longrightarrow}Z_q
\label{M_conv_dist_1}
\end{equation}
as $N\uparrow\infty$ if $q=\alpha$ and
\begin{equation}
N^{1-q/\alpha}\cdot\EMone\stackrel{d}{\longrightarrow}Z_q
\label{M_conv_dist_2}
\end{equation}
as $N\uparrow\infty$ if $q>\alpha$. Indeed, since $\{N_k/\tau\}_{k\ge
  1}$ is a subsequence of the integers for a given
$\tau\le\mathcal{T}$, combining (\ref{dist_Mt_M1}) with
(\ref{M_conv_dist_1}) and (\ref{M_conv_dist_2}) we get that
\begin{equation*}
\EMk-c\,\tau\ln\biggl(\frac{N_k}{\tau}\biggr)~\stackrel{d}{=}~\tau\cdot\Bigl[{}^q\!M^1_{N_k/\tau}-c\ln(N_k/\tau)\Bigr]\stackrel{d}{\longrightarrow}\tau\cdot Z_q
\end{equation*}
as $k\uparrow\infty$ if $q=\alpha$ and
\begin{equation*}
N_k^{1-q/\alpha}\cdot\EMk~\stackrel{d}{=}~\tau\cdot(N_k/\tau)^{1-q/\alpha}\cdot{}^q\!M^1_{N_k/\tau}\stackrel{d}{\longrightarrow}\tau\cdot Z_q
\end{equation*}
as $k\uparrow\infty$ if $q>\alpha$. The law of $Z_q$ is clearly
independent of $\tau$ as $\EMone$ does not depend on it.

Convergences in distribution (\ref{M_conv_dist_1}) and
(\ref{M_conv_dist_2}) follow from standard limit theory for sums of
i.i.d.\ random variables. We point out that $|X_1|^q$ belongs to the
normal domain of attraction of a stable variable with stability index
$\alpha/q\le 1$ when $q\ge\alpha$ (see \cite{Feller}, concluding
remark in Chapter XVII.5).  Indeed, $|X_1|^q$ is positive and
$\lim_{x\uparrow\infty}x^{\alpha/q}\cdot\prob[|X_1|^q>x]=c$ from
(\ref{tailabs}).  Then, since the increments $\{\Delta_n^1X\}_{n\ge
  0}$ over consecutive, non--overlapping time windows of size 1 form a
sequence of independent variables distributed as $X_1$ for L\'evy
processes, we realise that (see \cite{Feller}, Theorem 3 in Chapter
XVII.5)
\begin{equation}
\EMone-b_N=\frac{1}{N}\sum_{n=0}^{N-1}\bigl|\Delta_n^1X\bigr|^q-b_N\stackrel{d}{\longrightarrow}Z_q
\label{sumZa}
\end{equation}
as $N\uparrow\infty$ if $q=\alpha$ and
\begin{equation}
N^{1-q/\alpha}\cdot\EMone=\frac{1}{N^{\frac{q}{\alpha}}}\sum_{n=0}^{N-1}\bigl|\Delta_n^1X\bigr|^q\stackrel{d}{\longrightarrow}Z_q
\label{sumZq}
\end{equation}
as $N\uparrow\infty$ if $q>\alpha$, where $Z_q$ is a non--degenerate
stable variable of stability index $\alpha/q$ and $\{b_N\}_{N\ge 1}$
is any numerical sequence such that
\begin{equation*}
\lim_{N\uparrow\infty}\biggl\{\att\Bigl[N\sin\bigl(|X_1|^\alpha/N\bigr)\Bigr]-b_N\biggr\}
\end{equation*}
exists finite.  A centering procedure is required when $q=\alpha$ to
avoid a degenerate limit distribution.  The following lemma, whose
proof is reported in \ref{proof:lemma1}, states that $b_N$ can be
taken equal to $c\ln (N)$ for each $N\ge 1$, $c$ being the constant
defined by (\ref{c_def}). Theorem \ref{th:conv} is thus demonstrated.
\begin{lemma}
Let $c$ be the constant defined by
(\ref{c_def}). Then
\begin{equation*}
\lim_{N\uparrow\infty}\biggl\{\att\Bigl[N\sin\bigl(|X_1|^\alpha/N\bigr)\Bigr]-c\ln(N)\biggr\}
\end{equation*}
exists finite.
\end{lemma}

Proposition \ref{no_cost} is an immediate consequence of the following
important result due to Feller (see \cite{Durrett}, Theorem 2.5.9).
\begin{theorem}
\label{theorem_Feller}
Let $\{V_n\}_{n\ge 0}$ be a sequence of i.i.d.\ random variables
defined on some probability space $(\Omega,\mathcal{F},\prob)$ and let
$\{b_N\}_{N\ge 1}$ be an increasing sequence of positive numbers. If
$\att[|V_0|]=\infty$, then
\begin{equation*}
\lim_{N\uparrow\infty}\frac{1}{Nb_N}\sum_{n=0}^{N-1}V_n=0~~~~~~~\mbox{a.s.}
\end{equation*}
or 
\begin{equation*}
\limsup_{N\uparrow\infty}\frac{1}{Nb_N}\biggl|\sum_{n=0}^{N-1}V_n\biggr|=\infty~~~~~~~\mbox{a.s.}
\end{equation*}
according as the series $\sum_{N=1}^\infty\prob[|V_0|>N b_N]$ converges
or diverges.
\end{theorem}
In order to prove Proposition \ref{no_cost} we fix
$\tau\le\mathcal{T}$, we set $\lambda:=\lcm/\tau$, and we group the
first $N_k/\tau=\lambda k$ increments of the process over time windows
of size $\tau$ into $k$ consecutive blocks of $\lambda$ elements
each. This allows us to recast $\EMk$ as
\begin{equation*}
\EMk=\frac{1}{\lambda k}\sum_{n=0}^{\lambda k-1}\bigl|\Delta_{\tau n}^\tau X\bigr|^q
=\frac{1}{k}\sum_{i=0}^{k-1} V_i
\end{equation*}
with
\begin{equation*}
V_i:=\frac{1}{\lambda}\sum_{j=0}^{\lambda-1}\bigl|\Delta_{\tau(\lambda i+j)}^\tau X\bigr|^q.
\end{equation*}
Since $\{\Delta_{\tau n}^\tau X\}_{n\ge 0}$ is a sequence of
i.i.d.\ random variables, $\{V_i\}_{i\ge 0}$ is a sequence of
i.i.d.\ positive random variables. Moreover, as the $|\Delta_{\tau
  n}^\tau X|^q$'s are fat--tailed variables distributed as
$|X_\tau|^q$, for the sum of the first $\lambda$ of them we have that
(see \cite{Feller}, last corollary in Chapter VIII.8)
\begin{equation*}
\prob\Biggl[\sum_{j=0}^{\lambda-1}|\Delta_{\tau j}^\tau X|^q>x\Biggr]\sim \lambda\,\prob\bigl[|X_\tau |^q>x\bigr].
\end{equation*}
Thus, bearing in mind (\ref{tailabs}), we find that
\begin{equation*}
\prob[V_0>x]=\prob\Biggl[\sum_{j=0}^{\lambda-1}|\Delta_{\tau j}^\tau X|^q>\lambda x\Biggr]\sim \lambda\,\prob\bigl[|X_\tau|^q>\lambda x\bigr]
\sim \frac{\tau c\,\lambda^{1-\alpha/q}}{x^{\alpha/q}},
\end{equation*}
with the consequence that for each increasing sequence
$\{a_{N_k}\}_{k\ge 1}$ of positive numbers the series
$\sum_{k=1}^\infty\prob[V_0>k a_{N_k}]$ converges or diverges if and
only if the series $\sum_{k=1}^\infty(k a_{N_k})^{-\alpha/q}$
does. When $q\ge\alpha$, then $\att[V_0]=\infty$ and Theorem
\ref{theorem_Feller} applies stating that for a given increasing
sequence $\{a_{N_k}\}_{k\ge 1}$ of positive numbers
\begin{equation*}
\lim_{k\uparrow\infty}\frac{\EMk}{a_{N_k}}=\lim_{k\uparrow\infty}\frac{1}{k a_{N_k}}\sum_{i=0}^{k-1}V_i=0~~~~~~~\mbox{a.s.}
\end{equation*}
or 
\begin{equation*}
\limsup_{k\uparrow\infty}\frac{\EMk}{a_{N_k}}=\infty~~~~~~~\mbox{a.s.}
\end{equation*}
according as the series $\sum_{k=1}^\infty\prob[V_0>k a_{N_k}]$, and
hence $\sum_{k=1}^\infty(k a_{N_k})^{-\alpha/q}$, converges or
diverges. As $\tau$ is arbitrary and as the intersection of a finite
number of events with probability measure one has probability measure
one too, Proposition \ref{no_cost} is proven.

\subsection{Almost sure apparent multifractality}
\label{ourresult}

Numerical investigations carried out in \cite{CG2000} lead to the
conjecture that, irrespective of the order $q$, for almost all
$\omega\in\Omega$ the empirical moment $\EM(\omega)$ becomes
proportional to $\tau^{\nu_e(q)}$ at large $N$ with coefficient of
proportionality independent of $\tau$, $\nu_e(q)$ being the empirical
scaling function (\ref{empiricalnu}). Within this scenario the
coefficient of proportionality should be equal to $\EMone(\omega)$.
Then, we believe that a breakthrough in understanding apparent
multifractality of self--similar L\'evy processes may come from the
study of the ratio $\EM/\EMone$, with the purpose of demonstrating
that it converges a.s.\ to $\tau^{\nu_e(q)}$ as $N$ is sent to
infinite. Actually, this is an evident fact at orders $q$ smaller than
the stability index $\alpha$, where the strong law of large numbers
works. In contrast, it requires to hold in the regime $q\ge\alpha$
that a nontrivial decoupling between the randomness source $\omega$
and the size $\tau$ of time windows occurs at large $N$ as for almost all
$\omega$ the empirical moment $\EM(\omega)$ remains a fluctuating
quantity depending on the specific sample $\omega$, whose behaviour
cannot be regularised by the introduction of a norming constant as
stated by Proposition \ref{no_cost}.

In this paper we show that the conjecture holds for self--similar
L\'evy processes.  The reason of the decoupling between the randomness
source $\omega$ and the variable $\tau$ at orders $q\ge\alpha$ lies,
in a nutshell, in the fact that large values of the increments of the
process over time windows of integer size $\tau$ are due to large
fluctuations of only one of the $\tau$ corresponding unitary
increments. We let the following theorem to state our main result.
\begin{theorem}
\label{main}
Let $\Xcont$ be a self--similar L\'evy process of stability index
$\alpha<2$ defined on some probability space
$(\Omega,\mathcal{F},\prob)$ and let $\nu_e$ be the empirical scaling
function (\ref{empiricalnu}). Let $\mathcal{T}$ be a positive integer and let
{\upshape $N_k:=k\,\lcm$} be a multiple of {\upshape $\lcm$}. Then
\begin{equation*}
\prob\Biggl[\bigcap_{\tau=1}^{\mathcal{T}}\biggl\{\omega\in\Omega:\lim_{k\uparrow\infty}\frac{\EMk(\omega)}{\EMonek(\omega)}=\tau^{\nu_e(q)}\biggr\}\Biggr]=1
\end{equation*}
for all $q\ge 0$.
\end{theorem}

The theorem demonstrates that empirical moments possess their own
scaling law in time at any order $q$ for almost all the sample paths
of self--similar L\'evy processes. This fact leads to the systematic
emergence of apparent multifractality and completely explains the
numerical findings of \cite{SSL1999} and \cite{CG2000}. Apparent
multifractality emerges a.s.\ even when $q=\alpha$, although Theorem
\ref{th:conv} tells us that the distribution of estimators does not
scale with time in this case. We stress that the scaling law of
empirical moments needed to be disclosed of an appropriate stochastic
normalisation, which was able to balance fluctuations without hiding
the dependence on time. Theorem \ref{main} states that such a
stochastic normalisation is exactly $\EMone$.

Theorem \ref{main} is an immediate consequence of
(\ref{conv_emp_moment}) when $q<\alpha$. In order to prove the theorem
in the case $q\ge\alpha$ we observe that, given $N$ multiple of $\lcm$
and grouping the first $N$ unitary increments of the process into
$N/\tau$ consecutive blocks of $\tau\le\mathcal{T}$ elements each, the
estimator $\EMone$ can be rewritten as
\begin{equation*}
\EMone=\frac{1}{N}\sum_{n=0}^{N-1}\bigl|\Delta_n^1X\bigr|^q=\frac{1}{N}\sum_{n=0}^{N/\tau-1}\sum_{i=0}^{\tau-1}\bigl|\Delta_{\tau n+i}^1X\bigr|^q,
\end{equation*}
so that 
\begin{equation*}
\frac{\EM}{\EMone}=\tau \cdot \frac{\sum_{n=0}^{N/\tau-1}\bigl|\Delta_{\tau n}^\tau X\bigr|^q}{\sum_{n=0}^{N/\tau-1}\sum_{i=0}^{\tau-1}\bigl|\Delta_{\tau n+i}^1X\bigr|^q}.
\end{equation*}
Then, since $\{N_k/\tau\}_{k\ge 1}$ is a subsequence of the integers
for $\tau\le\mathcal{T}$ and since a finite intersection of measurable
sets with probability measure one has probability measure one too,
Theorem \ref{main} descends from the following proposition as
$\nu_e(q)=1$ for each $q\ge\alpha$.
\begin{proposition}
\label{lim_imp}
Assume that $\Xcont$ is a self--similar L\'evy process of stability
index $\alpha<2$ defined on some probability space
$(\Omega,\mathcal{F},\prob)$. Fix any $q\ge\alpha$ and integer $\tau\ge
1$. Then
\begin{equation*}
\lim_{N\uparrow\infty}\frac{\sum_{n=0}^{N-1}\bigl|\Delta_{\tau n}^\tau X\bigr|^q}{\sum_{n=0}^{N-1}\sum_{i=0}^{\tau-1}\bigl|\Delta_{\tau n+i}^1X\bigr|^q}=1~~~~~~~\mbox{a.s.}.
\end{equation*}
\end{proposition}

Proposition \ref{lim_imp} is interesting by its own and could be
easily extended to processes with more general stationary independent
fat--tailed increments. Its proof is outlined in the next
paragraph. This proposition states that the large fluctuations of
$|\Delta_{\tau n}^\tau X|^q$ are comparable to those of the sum
$\sum_{i=0}^{\tau-1}|\Delta_{\tau n+i}^1X|^q$ referring to the same
time windows, due to the hypothesis of independent increments. We
immediately recognise that this fact is possible only if a unique
unitary increment dominates by noticing that $\Delta_{\tau n}^\tau
X=\sum_{i=0}^{\tau-1}\Delta_{\tau n+i}^1X$. Thus, it is confirmed that
large values of the increments of the process over time windows of
size $\tau$ are ascribed to large fluctuations of only one of the
corresponding unitary increments if $q\ge\alpha$.

\subsection{Proof of Proposition \ref{lim_imp}}
\label{proofresult}

We outline the proof of Proposition \ref{lim_imp} postponing technical
details in the appendices with the purpose of not interrupting the
flow of the presentation.  We fix once and for all a moment order
$q\ge\alpha$ and, noticing that the instance $\tau=1$ is trivial, an
integer $\tau\ge 2$.  Since we found convenient reducing the proof to
a comparison among extreme events, at first we consider for each $n\ge
0$ the pointwise largest and second largest of the $\tau$ unitary
increments
$|\Delta_{\tau n}^1X|^q,|\Delta_{\tau n+1}^1X|^q,\ldots,|\Delta_{\tau n+\tau-1}^1X|^q$,
denoted by $U_n$ and $V_n$ respectively. Both $\{U_n\}_{n\ge 0}$ and
$\{V_n\}_{n\ge 0}$ are sequences of i.i.d.\ positive random variables
due the hypothesis of stationary independent increments. We observe
that $U_0$ is not larger than $x$ if and only if all the variables
$|\Delta_0^1X|^q,|\Delta_1^1X|^q,\ldots,|\Delta_{\tau-1}^1X|^q$ are
smaller than or equal to $x$, so that $\mathbb{P}[U_0\le
  x]=\prob^\tau[|\Delta_0^1X|^q\le x]=\prob^\tau[|X_1|^q\le x]$ and
$\mathbb{P}[U_0> x]\sim \tau\,\prob[|X_1|^q>x]$ follows. Similarly, $V_0$
is not larger than $x$ if and only if at least $\tau-1$ among these
variables are smaller than or equal to $x$, with the consequence that
$\mathbb{P}[V_0\le x]=\prob^\tau[|X_1|^q\le
  x]+\tau\,\prob[|X_1|^q>x]\,\prob^{\tau-1}[|X_1|^q\le x]$, providing
$\mathbb{P}[V_0>x]\sim (1/2)\tau(\tau-1)\,\prob^2[|X_1|^q>x]$. The
asymptotic formulas reflect a simple combinatorial argument: large
fluctuations of $U_0$ require large values of only one increment among
$|\Delta_0^1X|^q,|\Delta_1^1X|^q,\ldots,|\Delta_{\tau-1}^1X|^q$, while
large fluctuations of $V_0$ demand large values of at least two of
them. Combining these formulas with (\ref{tailabs}) we obtain that
\begin{equation}
\mathbb{P}[U_0>x]\sim\frac{\tau c}{x^\beta}
\label{Mtail}
\end{equation}
and
\begin{equation}
\mathbb{P}[V_0>x]\sim\frac{\tau(\tau-1)}{2}\frac{c^2}{x^{2\beta}},
\label{mtail}
\end{equation}
where $\beta:=\alpha/q\le 1$ and $c$ is given by (\ref{c_def}).  The
following lemma, which is proven in \ref{proof:lemma2}, brings the
$U_n$'s and $V_n$'s up showing that in order to prove Proposition
\ref{lim_imp} it is enough to demonstrate that the ratio
\begin{equation*}
R_N:=\frac{\sum_{n=0}^{N-1}V_n}{\sum_{n=0}^{N-1}U_n}
\end{equation*}
goes to zero a.s.\ as $N$ is sent to infinite. We notice that
$\sum_{n=0}^{N-1}U_n=0$ if and only if $\Delta_{n}^1X=0$ for each
$n\le \tau N-1$ and that both these events have probability measure equal
to zero since the distribution of $X_1$ does not put any finite mass
on any point.
\begin{lemma}
\label{lemma2}
There exists a positive constant $h$, in general depending on both $q$
and $\tau$, such that 
\begin{equation*}
\Biggl|\frac{\sum_{n=0}^{N-1}\bigl|\Delta_{\tau n}^\tau X\bigr|^q}{\sum_{n=0}^{N-1}\sum_{i=0}^{\tau-1}\bigl|\Delta_{\tau n+i}^1X\bigr|^q}-1\Biggr|
\le h\cdot (R_N)^e~~~~~~~\mbox{a.s.}
\end{equation*}
for all $N\ge 1$ with $e:=\min\{1,1/q\}$.
\end{lemma}

The fact that $\lim_{N\uparrow\infty}R_N=0$ a.s.\ is not surprising
from the point of view of intuition. As large fluctuations of the
$V_n$'s are much less probable than comparable fluctuations of the
$U_n$'s, one expects that the numerator of $R_N$ becomes smaller and
smaller in comparison with the denominator if $N$ becomes larger and
larger. However, one has to pay attention that there are subtle
samples for which the denominator of $R_N$ does not grow quickly
enough to dominate the numerator, with the consequence that $R_N$ is
not negligible in the large $N$ limit.  These samples have thus to be
isolated. Fortunately, they constitute a set which has a smaller and
smaller probability measure as $N$ is let to increase.  These
considerations prompt us to resort to the following strategy to prove
that $\lim_{N\uparrow\infty}R_N=0$ a.s.. For each $N\ge 1$ we
pick a positive number $\lambda_N$, which is for the moment
unknown. Then we consider the set $E$ of all $\omega\in\Omega$ with
the property that there exists $N_0\ge 1$ such that $\sum_{n=0}^{N-1}
U_n(\omega)>\lambda_N$ if $N\ge N_0$.  This set will correspond to the
set of samples for which the denominator of $R_N$ is able to dominate
the numerator. Finally we search for the $\lambda_N$'s that makes $E$
a set with probability measure equal to one on the one hand, and allow
to satisfy $\lim_{N\uparrow\infty}R_N(\omega)=0$ for almost all
$\omega\in E$ on the other hand.  In order to put into practice this
strategy, for each $N\ge 1$ we introduce the measurable set
\begin{equation*}
E_N:=\Biggl\{\omega\in\Omega:\sum_{n=0}^{N-1} U_n(\omega)>\lambda_N\Biggr\}.
\end{equation*}
The set $E$ of all $\omega\in\Omega$ for which there exists $N_0\ge 1$
with the property that $\sum_{n=0}^{N-1} U_n(\omega)>\lambda_N$ if
$N\ge N_0$ is nothing but the limit inferior of the sequence of the
$E_N$'s: $E=\bigcup_{n=1}^\infty\bigcap_{N=n}^\infty E_N$.  Indeed, if
$\omega\in E$ one can find $N_0\ge 1$ so that
$\omega\in\bigcap_{N=N_0}^\infty E_N$ and vice versa. We also
introduce the measurable set
\begin{equation*}
F:=\Biggl\{\omega\in\Omega:\lim_{N\uparrow\infty}\frac{1}{\lambda_N}\sum_{n=0}^{N-1} V_n(\omega)=0\Biggr\}
\end{equation*}
and observe that $\lim_{N\uparrow\infty}R_N(\omega)=0$ if $\omega\in
E\cap F$. Indeed, both the conditions $0\le
R_N(\omega)<(1/\lambda_N)\sum_{n=0}^{N-1}V_n(\omega)$ for all $N\ge
N_0$ with some $N_0\ge 1$ and
$\lim_{N\uparrow\infty}(1/\lambda_N)\sum_{n=0}^{N-1}V_n(\omega)=0$ are
satisfied when $\omega\in E\cap F$.  Thus, if a choice of the
$\lambda_N$'s for which both $\prob[E]=1$ and $\prob[F]=1$ exists,
then $\lim_{N\uparrow\infty}R_N=0$ a.s.\ follows.

A straight application of the {\it Borel--Cantelli Lemma} yields that
$\prob[E^c]=0$, and hence $\prob[E]=1$, if the condition
\begin{equation}
\sum_{N=1}^\infty\prob[E_N^c]=\sum_{N=1}^\infty\prob\Biggl[\sum_{n=0}^{N-1} U_n\le\lambda_N\Biggr]<\infty
\label{condition_l}
\end{equation}
is met. This amounts to say that the event where the denominator of
$R_N$ is not able to dominate the numerator is unlikely enough at
large $N$.  We let the following lemma, which is demonstrated in
\ref{proof:lemma3}, to introduce a choice of the $\lambda_N$'s that
fulfils such condition.
\begin{lemma}
\label{lemma3}
There exist two strictly positive constants $k$ and $K$, in general
depending on both $q$ and $\tau$, such that 
\begin{equation*}
\prob\Biggl[\sum_{n=0}^{N-1} U_n\le \lambda_N\Biggr]\le\frac{K}{N^2}
\end{equation*}
for all $N\ge 1$ with $\beta:=\alpha/q$ and
\begin{equation*}
\lambda_N:=\cases{
k N^{\frac{1}{\beta}}(\ln N+1)^{-\frac{1-\beta}{\beta}} & if $\beta<1$; \\
k N(\ln N+1) & if $\beta=1$.\\
}
\end{equation*}
\end{lemma}

The $\lambda_N$'s introduced by the lemma not only satisfy the
condition (\ref{condition_l}), in such a way that $\prob[E]=1$, but
even give that $\prob[F]=1$. Indeed, if the positive variable $V_0$
possesses finite expected value, namely if $q<2\alpha$ as stated by
(\ref{mtail}), then the strong law of the large numbers applies
ensuring that a measurable set $G$ exists with the properties that
$\prob[G]=1$ and
$\lim_{N\uparrow\infty}(1/N)\sum_{n=0}^{N-1}V_n(\omega)=\mathbb{E}[V_0]$
for all $\omega\in G$. It follows that
$\lim_{N\uparrow\infty}(1/\lambda_N)\sum_{n=0}^{N-1}V_n(\omega)=0$ for
each $\omega\in G$ since $N/\lambda_N$ approaches zero as $N$ is sent
to infinite. Thus $G\subseteq F$ and $\prob[F]=1$. If instead $V_0$
does not possess finite expected value, namely if $q\ge 2\alpha$, then
Theorem \ref{theorem_Feller} applies yielding that
$\lim_{N\uparrow\infty}(1/\lambda_N)\sum_{n=0}^{N-1}V_n=0$ a.s.\ since
one can easily verify that $\{\lambda_N/N\}_{N\ge 1}$ is an increasing
sequence of positive numbers and that
$\sum_{N=1}^\infty\lambda_N^{-2\beta}<\infty$. The latter fact amounts
to $\sum_{N=1}^\infty\prob[V_0>\lambda_N]<\infty$ thanks to
(\ref{mtail}).  The proof of Proposition \ref{lim_imp} is thus
concluded.

\section{Conclusions}
\label{conclusion}

In this paper we have proven that empirical moments of self--similar
L\'evy processes possess almost surely their own scaling laws with a
piecewise--linear scaling function. Such scaling laws differ from
those of the underlying process, as they extend the range of validity
of the latter that is limited because of the presence of fat tails.
Since only the scaling properties of empirical moments can be detected
from a particular sample, our results explain the emergence of
apparent multifractality in self--similar L\'evy processes once and
for all.

The main lesson to be learnt from this study is that inferring
multifractality from scaling laws of empirical moments may lead to
mistaken conclusions in many practical situations where fat--tailed
variables are involved.  Nonlinearities of the scaling function should
not be taken as an evidence for genuine multifractality unless the
range of existing moments is known, otherwise better tests for
multiplicative noise should be developed to assess whether the data
really account for such a property or not.

Results proven in this paper are concerned with processes with
independent increments. As this condition is usually not met in real
time series, dependence should be considered in future research in
order to clarify the interplay between multifractality and
heavy--tails for more general stochastic processes.  The quest is
particularly important for financial data, which often display fat
tails and strong dependence, and which produce nonlinear scaling
functions at the same time. We expect that the use of a stochastic
normalisation to regularise empirical moments, represented by the
estimator of moments at a given size of time windows, can help to
explain the emergence of apparent multifractality, if any, in a wide
context. The particular case of self--similar L\'evy processes
constitutes an example in which such procedure succeeds.

\section*{Acknowledgement}
The author is grateful to Fulvio Baldovin and Attilio Stella for
stimulating discussions which led to the initiation of this study.

\appendix
\section*{Appendices}

\section{Proof of Lemma 1}
\label{proof:lemma1}

Applying the {\it Fubini--Tonelli Theorem} to the integral
\begin{equation*}
\int_0^\infty dx\int_{\Omega}\prob[d\omega]~\cos(x/N)\mathbbm{1}\bigl(|X_1(\omega)|^\alpha>x\bigr)
\end{equation*}
we get that
\begin{equation*}
\att\Bigl[N\sin\bigl(|X_1|^\alpha/N\bigr)\Bigr]=\int_0^\infty dx~\cos(x/N)\prob\bigl[|X_1|^\alpha>x\bigr].
\end{equation*}
Then, subtracting and adding the convergent {\it Dirichlet--type
  integral} $\int_1^\infty dx~\frac{\cos(x/N)}{x}$ we can write
\begin{eqnarray}
\nonumber
\att\Bigl[N\sin\bigl(|X_1|^\alpha/N\bigr)\Bigr]-c\ln(N)&=&\int_0^1 dx~\cos(x/N)\prob\bigl[|X_1|^\alpha>x\bigr]\\
\nonumber
&+&\int_1^\infty dx~\cos(x/N)\prob\bigl[|X_1|^\alpha>x\bigr]-c\ln(N)\\
\nonumber
&=&\int_0^1 dx~\cos(x/N)\prob\bigl[|X_1|^\alpha>x\bigr]\\
\nonumber
&+&\int_1^\infty dx~\cos(x/N)\biggl\{\prob\bigl[|X_1|^\alpha>x\bigr]-\frac{c}{x}\biggr\}\\
&+&c\Biggl\{\int_1^\infty dx~\frac{\cos(x/N)}{x}-\ln(N)\Biggr\}.
\label{lemma11}
\end{eqnarray}
Since we can find $d>0$ in such a way that
$|\prob[|X_1|>x]-c/x^\alpha|\le d/x^{2\alpha}$ for all $x\ge 1$ when
$c$ is as in (\ref{c_def}) (see \cite{Zolotarev}, Theorem 2.4.2 if
$\alpha<1$ and Corollary 2 of Theorem 2.5.1 if $\alpha>1$), the {\it
  Dominated Convergence Theorem} yields that the first and the second
terms of (\ref{lemma11}) form convergent sequences. Thus, the trivial
identity
\begin{equation*}
\int_1^\infty dx~\frac{\cos(x/N)}{x}-\ln(N)=\int_{1/N}^1 dx~\frac{\cos(x)-1}{x}+\int_1^\infty dx~\frac{\cos(x)}{x}
\end{equation*}
allows us to conclude that
\begin{eqnarray}
\nonumber
\lim_{N\uparrow\infty}\biggl\{\att\Bigl[N\sin\bigl(|X_1|^\alpha/N\bigr)\Bigr]-c\ln(N)\biggr\}
&=&\int_0^1 dx~\prob\bigl[|X_1|^\alpha>x\bigr]\\
\nonumber
&+&\int_1^\infty dx~\biggl\{\prob\bigl[|X_1|^\alpha>x\bigr]-\frac{c}{x}\biggr\}\\
\nonumber
&+&c\Biggl\{\int_0^1 dx~\frac{\cos(x)-1}{x}+\int_1^\infty dx~\frac{\cos(x)}{x}\Biggr\}.
\end{eqnarray}

\section{Proof of Lemma 2}
\label{proof:lemma2}

The proof of Lemma \ref{lemma2} relies on the following bounds, which
hold for all numbers $0\le v\le u$ and $\tau\ge 1$ with $h:=2\tau$ or
$h:=(q+2)\tau^q$ according as $q\le 1$ or $q>1$:
\begin{equation}
\Bigl[\max\Bigl\{0,u^{\frac{1}{q}}-(\tau-1)v^{\frac{1}{q}}\Bigr\}\Bigr]^q \ge \cases{
u+(\tau-1)v-h\,v & if $q\le 1$; \\
u+(\tau-1)v-h\,u^{1-\frac{1}{q}}\,v^{\frac{1}{q}} & if $q>1$\\
}~
\label{lemma21}
\end{equation}
and
\begin{equation}
\Bigl[u^{\frac{1}{q}}+(\tau-1)v^{\frac{1}{q}}\Bigr]^q \le \cases{
u+h\,v & if $q\le 1$; \\
u+h\,u^{1-\frac{1}{q}}\,v^{\frac{1}{q}} & if $q>1$.\\
}
\label{lemma22}
\end{equation}
We shall prove these inequalities later.  At the moment, we exploit
them to show that for the given order $q>0$ and integer $\tau\ge 2$ and
for each $N\ge 1$ and real numbers
$\delta_0,\delta_1,\ldots,\delta_{\tau N-1}$ not all equal to zero
\begin{equation}
1-h\cdot\Biggl(\frac{\sum_{n=0}^{N-1}v_n}{\sum_{n=0}^{N-1}u_n}\Biggr)^e\le 
\frac{\sum_{n=0}^{N-1}\bigl|\sum_{i=0}^{\tau-1}\delta_{\tau n+i}\bigr|^q}{\sum_{n=0}^{N-1}\sum_{i=0}^{\tau-1}|\delta_{\tau n+i}|^q} \le 
1+h\cdot\Biggl(\frac{\sum_{n=0}^{N-1}v_n}{\sum_{n=0}^{N-1}u_n}\Biggr)^e,~~~
\label{lemma2main}
\end{equation}
where $u_n$ and $v_n$ are used to denote the largest and the second
largest of the $\tau$ values $|\delta_{\tau n}|^q,|\delta_{\tau
  n+1}|^q,\ldots,|\delta_{\tau n+\tau-1}|^q$ respectively and
$e:=\min\{1,1/q\}$.  The lemma follows from this formula setting
$\delta_n:=\Delta_n^1X(\omega)$, with $\omega\in\Omega$ chosen so that
the $\Delta_n^1X(\omega)$'s are not all equal to zero, and recognising
that $\sum_{i=0}^{\tau-1}\delta_{\tau
  n+i}=\sum_{i=0}^{\tau-1}\Delta_{\tau n+i}^1X(\omega)=\Delta_{\tau
  n}^\tau X(\omega)$.

In order to prove (\ref{lemma2main}), we take advantage of the fact
that $u_n^{1/q}$ and $v_n^{1/q}$ are the first and the second largest
of $|\delta_{\tau n}|,|\delta_{\tau n+1}|,\ldots,|\delta_{\tau
  n+\tau-1}|$ to obtain the straight bounds
\begin{equation*}
\max\Bigl\{0,u_n^{\frac{1}{q}}-(\tau-1)v_n^{\frac{1}{q}}\Bigr\}\le\Biggl|\sum_{i=0}^{\tau-1}\delta_{\tau n+i}\Biggr|\le u_n^{\frac{1}{q}}+(\tau-1)v_n^{\frac{1}{q}}.
\end{equation*}
Then, as $0\le v_n\le u_n$, we invoke (\ref{lemma21}) and (\ref{lemma22}) with $u:=u_n$
and $v:=v_n$ to state that for each $n<N$ on the one hand
\begin{equation}
\Biggl|\sum_{i=0}^{\tau-1}\delta_{\tau n+i}\Biggr|^q \ge\cases{
u_n+(\tau-1)v_n-h\,v_n & if $q\le 1$; \\
u_n+(\tau-1)v_n-h\,u_n^{1-\frac{1}{q}}\,v_n^{\frac{1}{q}} & if $q>1$,\\
}
\label{lemma23}
\end{equation}
and on the other hand
\begin{equation}
\Biggl|\sum_{i=0}^{\tau-1}\delta_{\tau n+i}\Biggr|^q \le\cases{
u_n+h\,v_n & if $q\le 1$; \\
u_n+h\,u_n^{1-\frac{1}{q}}\,v_n^{\frac{1}{q}} & if $q>1$.\\
}
\label{lemma24}
\end{equation}
At this point, carrying out the sum over $n$ is what remains to be
done.  In this respect, we observe that $1-1/q>0$ when $q>1$ so that
the {\it H\"{o}lder's inequality} applies ensuring that
\begin{equation*}
\sum_{n=0}^{N-1}u_n^{1-\frac{1}{q}}\,v_n^{\frac{1}{q}}\le\Biggl(\sum_{n=0}^{N-1}u_n\Biggr)^{1-\frac{1}{q}}\Biggl(\sum_{n=0}^{N-1}v_n\Biggr)^{\frac{1}{q}}=
\sum_{n=0}^{N-1}u_n\,\Biggl(\frac{\sum_{n=0}^{N-1}v_n}{\sum_{n=0}^{N-1}u_n}\Biggr)^{\frac{1}{q}}.
\end{equation*}
We notice that $\sum_{n=0}^{N-1}u_n$ is different from zero due to the
hypothesis that the $\delta_n$'s are not all equal to zero.  Carrying
out the sum over $n$ in (\ref{lemma23}) and (\ref{lemma24}) and
setting $S_N:=\sum_{n=0}^{N-1}[u_n+(\tau-1)v_n]$ and
$T_N:=\sum_{n=0}^{N-1}u_n$ we reach the result
\begin{equation*}
S_N-h\,T_N\Biggl(\frac{\sum_{n=0}^{N-1}v_n}{\sum_{n=0}^{N-1}u_n}\Biggr)^e\le
\sum_{n=0}^{N-1}\Biggl|\sum_{i=0}^{\tau-1}\delta_{\tau n+i}\Biggr|^q
\le T_N\Biggl[1+h\Biggl(\frac{\sum_{n=0}^{N-1}v_n}{\sum_{n=0}^{N-1}u_n}\Biggr)^e\Biggr],
\end{equation*}
where $e:=\min\{1,1/q\}$. This result proves (\ref{lemma2main}) once
combined with the fact that
\begin{equation*}
T_N\le\sum_{n=0}^{N-1}\sum_{i=0}^{\tau-1}|\delta_{\tau n+i}|^q\le S_N
\end{equation*}
as $u_n\le\sum_{i=0}^{\tau-1}|\delta_{\tau n+i}|^q\le u_n+(\tau-1)v_n$ for all
$n$ by the definition of $u_n$ and $v_n$.

We conclude the proof of the lemma showing the validity of
(\ref{lemma21}) and (\ref{lemma22}). The case $u=0$ is trivial since
$v=0$ as $0\le v\le u$. The case $u>0$ descends from the inequalities
\begin{equation}
1-\xi^q\le \bigl(\max\{0,1-\xi\}\bigr)^q\le (1+\xi)^q\le 1+\xi^q
\label{lemma26}
\end{equation}
if $q\le 1$ and
\begin{equation}
1-\xi^q-q\xi\le \bigl(\max\{0,1-\xi\}\bigr)^q\le (1+\xi)^q\le 1+q (1+\xi)^{q-1}\xi
\label{lemma27}
\end{equation}
when $q>1$, which hold for each positive number $\xi$. Verifying such
inequalities is a simple exercise of calculus, so we omit the details
and move on. The instance $q\le 1$ of (\ref{lemma21}) and
(\ref{lemma22}) follows from (\ref{lemma26}), which after setting
$\xi:=(\tau-1)(v/u)^{1/q}$ and multiplying by $u$ yields
\begin{eqnarray}
\nonumber
\Bigl[\max\Bigl\{0,u^{\frac{1}{q}}-(\tau-1)v^{\frac{1}{q}}\Bigr\}\Bigr]^q &\ge& u-(\tau-1)^q v\\
\nonumber
&=&u+(\tau-1)v-\bigl[\tau-1+(\tau-1)^q\bigr]v\\
\nonumber
&\ge& u+(\tau-1)v-h\,v
\end{eqnarray}
and 
\begin{equation*}
\Bigl[u^{\frac{1}{q}}+(\tau-1)v^{\frac{1}{q}}\Bigr]^q \le  u+(\tau-1)^q v \le u+h\,v
\end{equation*}
since both $\tau-1+(\tau-1)^q\le 2\tau$ and $(\tau-1)^q\le 2\tau$ as
$q\le 1$. The instance $q>1$ of (\ref{lemma21}) and (\ref{lemma22})
follows from (\ref{lemma27}) and the hypothesis $v\le u$, which in
particular entails that $v/u\le (v/u)^{1/q}$ and $v\le
u^{1-1/q}v^{1/q}$ as a consequence since $v/u\le 1$ and
$1/q<1$. Setting $\xi:=(\tau-1)(v/u)^{1/q}$ in (\ref{lemma27}),
multiplying by $u$, and bearing in mind that $v\le u^{1-1/q}v^{1/q}$
and that $v/u\le 1$ we find that
\begin{eqnarray}
\nonumber
\Bigl[\max\Bigl\{0,u^{\frac{1}{q}}-(\tau-1)v^{\frac{1}{q}}\Bigr\}\Bigr]^q &\ge& u-(\tau-1)^q v-q(\tau-1)u^{1-\frac{1}{q}}v^{\frac{1}{q}}\\
\nonumber
&=&u+(\tau-1)v-\bigl[\tau-1+(\tau-1)^q\bigr]v-q(\tau-1)u^{1-\frac{1}{q}}v^{\frac{1}{q}}\\
\nonumber
&\ge&u+(\tau-1)v-\bigl[(q+1)(\tau-1)+(\tau-1)^q\bigr]u^{1-\frac{1}{q}}v^{\frac{1}{q}}\\
\nonumber
&\ge& u+(\tau-1)v-h\,u^{1-\frac{1}{q}}\,v^{\frac{1}{q}}
\end{eqnarray}
and
\begin{eqnarray}
\nonumber
\Bigl[u^{\frac{1}{q}}+(\tau-1)v^{\frac{1}{q}}\Bigr]^q &\le & u+q(\tau-1)\biggl[1+(\tau-1)\Bigl(\frac{v}{u}\Bigr)^{\frac{1}{q}}\biggr]^{q-1} u^{1-\frac{1}{q}}v^{\frac{1}{q}}\\
\nonumber
&\le& u+q(\tau-1)\tau^{q-1} u^{1-\frac{1}{q}}v^{\frac{1}{q}}\le u+h\, u^{1-\frac{1}{q}}\,v^{\frac{1}{q}}
\end{eqnarray}
since both $(q+1)(\tau-1)+(\tau-1)^q\le (q+2)\tau^q$ and $q(\tau-1)\tau^{q-1}\le
(q+2)\tau^q$ as $q>1$.

\section{Proof of Lemma 3}
\label{proof:lemma3}

Given two numbers $\lambda>0$ and $\xi>0$ the {\it Markov's
  inequality} yields 
\begin{eqnarray}
\nonumber
\prob\Biggl[\sum_{n=0}^{N-1} U_n\le \lambda \Biggr]&=&
\prob\Biggl[\exp\Biggl(-\xi\sum_{n=0}^{N-1} U_n\Biggr)\ge \exp(-\xi\lambda) \Biggr]\\
\nonumber
&\le& \mbox{e}^{\xi\lambda}\cdot\att\Biggl[\exp\Biggl(-\xi\sum_{n=0}^{N-1} U_n\Biggr)\Biggr]=\mbox{e}^{\xi\lambda}\cdot\att^N\bigl[\exp(-\xi U_0)\bigr],
\end{eqnarray}
where the fact that the $U_n$'s are i.i.d.\ variables has been used to
obtain the last equality. The proof of the lemma moves from this
inequality and is based on the property of $U_0$ that there exist two
constants $\delta>0$ and $\eta>0$, in general depending on both
$\beta$ and $\tau$, such that the bound
\begin{equation}
\att\bigl[\exp(-\xi U_0)\bigr]\le\cases{
\mbox{e}^{-\delta \xi^\beta} & if $\beta<1$; \\
\mbox{e}^{\delta \xi\ln \xi} & if $\beta=1$\\
}
\label{expbound:1}
\end{equation}
holds for positive $\xi\le\eta$. This bound follows from (\ref{Mtail})
and will be proven later. Now we observe that combining the Markov's
inequality with such bound we get that if $0<\xi\le\eta$, then
\begin{equation}
\prob\Biggl[\sum_{n=0}^{N-1} U_n\le \lambda \Biggr]\le\cases{
\mbox{e}^{\xi\lambda-\delta N \xi^\beta} & if $\beta<1$; \\
\mbox{e}^{\xi\lambda +\delta N \xi\ln \xi} & if $\beta=1$.\\
}
\label{expbound:2}
\end{equation}

The lemma is a straight consequence of (\ref{expbound:2}) if the right
choice of $\xi$ and $\lambda$ is made. For each $N\ge 1$ we set
\begin{equation*}
\xi_N:=\cases{
\Bigl[\frac{2}{\delta (1-\beta)}\Bigr]^{\frac{1}{\beta}} N^{-\frac{1}{\beta}}(\ln N+1)^{\frac{1}{\beta}} & if $\beta<1$; \\
(\mbox{e}N^2)^{-\frac{1}{3}} & if $\beta=1$\\
}
\end{equation*}
and we notice that since $\lim_{N\uparrow\infty}\xi_N=0$ we can find
$N_0\ge 1$ with the property that $\xi_N\le\eta$ and
$\delta(N/\mbox{e})^{\frac{1}{3}}\ge 6$ for all $N\ge N_0$. Then,
taking $\lambda_N$ as in the statement of the lemma with
\begin{equation*}
k:=\cases{
\beta\, \delta^{\frac{1}{\beta}} \Bigl(\frac{1-\beta}{2}\Bigr)^{\frac{1-\beta}{\beta}} & if $\beta<1$; \\
\frac{\delta}{3} & if $\beta=1$,\\
}
\end{equation*}
simple algebra shows that $\xi_N\lambda_N-\delta N \xi_N^\beta=-2(\ln
N+1)\le -2\ln N$ when $\beta<1$ and that $\xi_N\lambda_N +\delta N
\xi_N\ln \xi_N=-(1/3)\delta(N/\mbox{e})^{\frac{1}{3}}\ln N\le -2\ln N$
if $\beta=1$ and $N\ge N_0$. Thus, plugging $\lambda:=\lambda_N$ and
$\xi:=\xi_N$ in (\ref{expbound:2}) we reach the result
\begin{equation*}
\prob\Biggl[\sum_{n=1}^N U_n\le \lambda_N \Biggr]\le\frac{1}{N^2}
\end{equation*}
for all $N\ge N_0$ and the lemma is proven with $K:=N_0^2$.

To conclude, we show the validity of the bound (\ref{expbound:1}).  To
begin with, given $\xi>0$, we apply the Fubini--Tonelli Theorem to the
integral
\begin{equation*}
\int_0^\infty dx\int_{\Omega}\prob[d\omega]~\exp(-\xi x)\mathbbm{1}\bigl(U_0(\omega)>x\bigr)
\end{equation*}
in order to obtain that
\begin{equation*}
\att\bigl[\exp(-\xi U_0)\bigr]=1-\xi\int_0^\infty dx~\exp(-\xi x)\prob\bigl[U_0>x\bigr].
\end{equation*}
Then, as (\ref{Mtail}) entails that $x^\beta\cdot\prob[U_0>x]\ge
(1/2)\tau c$ for all $x\ge s$ with some $s\ge 0$, in general depending
on both $\beta$ and $\tau$, we find that
\begin{eqnarray}
\nonumber
\att\bigl[\exp(-\xi U_0)\bigr]&\le&1-\xi\int_s^\infty dx~\exp(-\xi x)\prob\bigl[U_0>x\bigr]\\
&\le&1-\frac{\tau c\,\xi^{\beta}}{2}\int_{\xi s}^\infty y^{-\beta}\mbox{e}^{-y}dy,
\label{expbound:3}
\end{eqnarray}
where a change of variable has been performed.  On the other hand,
there exists a positive number $\eta$, in general depending on both
$\beta$ and $\tau$, with the property that for $\xi\le \eta$ one has
$\int_{\xi s}^\infty y^{-\beta}\mbox{e}^{-y}dy\ge \Gamma(1-\beta)/2$
if $\beta<1$, $\Gamma$ being the Euler's gamma function, and
$\int_{\xi s}^\infty y^{-\beta}\mbox{e}^{-y}dy\ge -\ln \xi/2$ when
$\beta=1$ since $\lim_{\xi\downarrow 0}\int_{\xi s}^\infty
y^{-\beta}\mbox{e}^{-y}dy=\Gamma(1-\beta)$ if $\beta<1$ and
$\lim_{\xi\downarrow 0}-(1/\ln\xi)\int_{\xi s}^\infty
y^{-1}\mbox{e}^{-y}dy=1$.  Thus, combining (\ref{expbound:3}) with
these bounds we find that for each positive $\xi\le\eta$
\begin{equation*}
\att\bigl[\exp(-\xi U_0)\bigr]\le\cases{
1-\delta \xi^\beta & if $\beta<1$; \\
1+\delta \xi\ln \xi & if $\beta=1$,\\
}
\end{equation*}
where $\delta:=\tau c\,\Gamma(1-\beta)/4>0$ if $\beta<1$ and
$\delta:=\tau c/4>0$ when $\beta=1$.  This result gives (\ref{expbound:1})
as $1+\zeta\le\mbox{e}^\zeta$ for any $\zeta$.

\end{document}